\pdfoutput=1
\documentclass{amsart}
\usepackage{amsmath}
\usepackage{amssymb}
\usepackage[mathscr]{euscript}
\usepackage{rotating}
\usepackage[all]{xy}
\usepackage{hyperref}

\newtheorem{theorem}{Theorem}[section]
\newtheorem{prop}[theorem]{Proposition}

\theoremstyle{definition}

\newtheorem{definition}{Definition}
\newtheorem{remark}{Remark}

\newcommand\C{\mathbb{C}}

\newcommand\R{\mathbb{R}}
\newcommand\Z{\mathbb{Z}}
\newcommand\Q{\mathbb {Q}}
\newcommand\Fi{\mathbb{F}}
\newcommand\Af{\mathbb{A}_f}
\newcommand\Ade{\mathbb{A}}

\newcommand\Hp{{}^{p}\mathrm{H}}
\newcommand\Ho{\mathrm{H}}
\newcommand\IH{\mathrm{IH}}
\newcommand\IC{{\bf IC}}

\newcommand\G{{\bf G}}
\newcommand\GU{{\bf GU}}
\newcommand\GSp{{\bf GSp}}

\newcommand\GL{{\bf GL}}

\newcommand\Ar{\mathrm{A}}
\newcommand\K{\mathrm {K}}

\newcommand\A{{\bf A}}
\newcommand\M{{\bf M}}

\newcommand\F{{\mathcal F}}

\newcommand\ggoth{\mathfrak{g}}

\newcommand\Scal{\mathcal{S}}
\newcommand\pp{\mathfrak{p}}
\newcommand\Of{\mathcal{O}}
\DeclareMathOperator{\Spec}{Spec}
\DeclareMathOperator{\Ima}{Im}
\DeclareMathOperator{\Gal}{Gal}
\newcommand\Hecke{\mathcal{H}}

\newcommand\DP{{}^{w}D}

\newcommand\sous{\setminus}
\newcommand\til{\widetilde}
\newcommand\X{{\mathcal X}}

\newcommand\fl{\longrightarrow}
\newcommand\fle{\longmapsto}
\newcommand\iso{\stackrel {\sim} {\fl}}
\newcommand\limik {\displaystyle{\lim_{\overrightarrow{\scriptstyle \,\,\,\K\,\,\,}}}}
\newcommand\ungras{1\!\!\mkern -1mu1}

\begin{document}

\title[The intersection complex as a weight truncation]{The intersection complex as a weight truncation and an application
to Shimura varieties}
\author[]{Sophie Morel}
\thanks{This text was written while I was working as a Professor
at the Harvard mathematics
department and supported by the Clay Mathematics Institute as a Clay
Research Fellow. I would like to thank the referee for their useful comments
about the first version of this paper.}

\address{Department of Mathematics, Harvard University,
One Oxford Street, Cambridge, MA 02138, USA}
\curraddr{Princeton Mathematics Department, Fine Hall, Washington Road, Princeton NJ 08544-1000, USA}
\email{smorel@math.princeton.edu}

\begin{abstract}
The purpose of this talk is to present an (apparently) new way to look at the
intersection complex of a singular variety over a finite field,
or, more generally, at the
intermediate extension functor on pure perverse sheaves, and an application
of this to the cohomology of noncompact Shimura varieties.

\noindent MSC (2000) : Primary 11F75; Secondary 11G18, 14F20.

\end{abstract}

\maketitle

\tableofcontents


\section{Shimura varieties}
\label{section:SV}

\subsection{The complex points}
\label{subsection:SV1}

In their simplest form, Shimura varieties are just locally symmetric varieties
associated to certain connected reductive groups over $\Q$. So let
$\G$ be a connected reductive group over $\Q$ satisfying the conditions in
1.5 of Deligne's article \cite{D}. To be precise, we are actually fixing $\G$
and a morphism $h:\C^\times\fl\G(\R)$ that is algebraic over $\R$.
Let us just remark here that these conditions
are quite restrictive.
For example, they exclude the group
$\GL_n$ as soon as $n\geq 3$. The groups $\G$ that we want to think about
are, for example,
the group $\GSp_{2n}$ (the general symplectic group of a symplectic space of
dimension $2n$ over $\Q$) or the general unitary group of a hermitian space
over a quadratic imaginary extension of $\Q$. The conditions on $\G$ ensure
that the symmetric space $\X$ of $\G(\R)$ is a hermitian symmetric domain;
so $\X$ has a canonical complex structure.
Remember that $\X=\G(\R)/\K'_\infty$, where $\K'_\infty$ is the centralizer in
$\G(\R)$ of $h(\C^\times)$. In the examples we consider, $\K'_\infty$ is the
product of a maximal compact subgroup $\K_\infty$
of $\G(\R)$ and of $\Ar_\infty:=\A(\R)^0$,
where $\A$ is the maximal $\Q$-split torus of the center of $\G$. (To avoid
technicalities, many authors assume that the maximal $\R$-split torus in the
center of $\G$ is also $\Q$-split. We will do so too.)

The locally symmetric spaces associated to $\G$ are the quotients
$\Gamma\sous\G(\R)$, where $\Gamma$ is an \emph{arithmetic subgroup}
of $\G(\Q)$,
that is, a subgroup of $\G(\Q)$ such that, for some (or any) $\Z$-structure
on $\G$, $\Gamma\cap\G(\Z)$ is of finite index in $\Gamma$ and
in $\G(\Z)$. If $\Gamma$ is small enough (for example, if it is torsion-free),
then $\Gamma\sous\X$ is a smooth
complex analytic variety. In fact, by the work of Baily and Borel
(\cite{BB}), it is even a quasi-projective algebraic variety.

In this text, we prefer to use the adelic point of view, as it leads to 
somewhat simpler
statements. So let $\K$ be a compact open subgroup of $\G(\Af)$, where
$\Af=\widehat{\Z}\otimes_\Z\Q$ is the ring of finite adeles of $\Q$.
This means that $\K$ is a subgroup of $\G(\Af)$ such that, for some (or any)
$\Z$-structure on $\G$, $\K\cap\G(\widehat{\Z})$ is of finite index in
$\K$ and in $\G(\widehat{\Z})$. Set
\[S^\K(\C)=\G(\Q)\sous(\X\times\G(\Af)/\K),\]
where $\G(\Q)$ acts on $\X\times\G(\Af)/\K$ by the formula $(\gamma,(x,g\K))
\fle (\gamma\cdot x,\gamma g\K)$.

This space $S^\K(\C)$ is related to the previous quotients $\Gamma\sous\X$ in
the following way. By the strong approximation theorem, $\G(\Q)\sous\G(\Af)/\K$
is finite. Let $(g_i)_{i\in I}$ be a finite family in $\G(\Af)$ such that
$\G(\Af)=\coprod_{i\in I}\G(\Q)g_i\K$. For every $i\in I$, set
$\Gamma_i=\G(\Q)\cap g_i\K g_i^{-1}$. Then the $\Gamma_i$ are arithmetic
subgroups of $\G(\Q)$, and
\[S^\K(\C)=\coprod_{i\in I}\Gamma_i\sous\X.\]

In particular, we see that, if $\K$ is small enough, then $S^\K(\C)$ is the
set of complex points of a smooth quasi-projective complex algebraic variety,
that we will denote by $S^\K$. These are the \emph{Shimura varieties}
associated to $\G$ and $h:C^\times\fl\G(\R)$ (over $\C$).
From now on, we will assume always that the group $\K$ is
small enough.

\begin{remark} If $\G=\GL_2$, then $S^\K$ is a modular curve, or rather, a
finite disjoint union of modular curves; it parametrizes elliptic curves
with a certain level structure (depending on $\K$).
Higher-dimensional generalizations of this
are the Shimura varieties for the symplectic groups $\G=\GSp_{2n}$; they
are called the Siegel modular varieties, and parametrize principally
polarized abelian varieties with a level structure (depending on $\K$).
Some other Shimura varieties have been given a name.
For example, if $\G$ is the general unitary group of a $3$-dimensional
hermitian vector space $V$ over an imaginary quadratic extension of
$\Q$ such that $V$ has signature $(2,1)$ at infinity, then $S^\K$ is called
a Picard modular surface.

\end{remark}

\subsection{The projective system and Hecke operators}

If $\K'\subset\K$ are two open compact subgroups of $\G(\Af)$, then there
is an obvious projection $S^{\K'}(\C)\fl S^\K(\C)$, and it defines a finite
{\'e}tale morphism $S^{\K'}\fl S^{\K}$; if $\K'$ is normal in $\K$, then this
morphism is Galois, with Galois group $\K/\K'$. So we can see the Shimura
varieties $S^\K$ as a projective system $(S^\K)_{\K\subset\G(\Af)}$ indexed by
(small enough) open compact subgroups of $\G(\Af)$, and admitting a right
continuous action of $\G(\Af)$.

More generally, if $\K',\K$ are two open compact subgroups of $\G(\Af)$ and
$g\in\G(\Af)$, then we get a correspondence $[\K'g\K]:S^{\K\cap g^{-1}\K'g}
\fl S^\K\times S^{\K'}$ in the following way. The first map is the obvious
projection $S^{\K\cap g^{-1}\K'g}\fl S^\K$, and the second map is
the composition
of the obvious projection $S^{\K\cap g^{-1}\K'g}\fl S^{g^{-1}\K'g}$ and of the
isomorphism $S^{g^{-1}\K'g}\iso S^{\K'}$. This is the
\emph{Hecke correspondence}
associated to $g$ (and $\K,\K'$).

Let $\Ho^*$ be a cohomology theory with coefficients in a ring $A$
that has good fonctoriality properties (for
example, Betti cohomology with coefficients in $A$) and $\K$
be an open compact subgroup of $\G(\Af)$. Then
the Hecke correspondences define an action of the Hecke algebra at level $\K$,
$\Hecke_\K(A):=C(\K\sous\G(\Af)/\K,A)$ (of bi-$\K$-invariant functions from
$\G(\Af)$ to $A$, with the algebra structure given
by the convolution product), on the cohomology $\Ho^*(S^\K)$. For every
$g\in\G(\Af)$, we make
$\ungras_{\K g\K}\in\Hecke_\K(A)$ act by the correspondence $[\K g^{-1}\K]$.

Let $\Hecke(A)=\bigcup_\K\Hecke_\K(A)=C_c^\infty(\G(\Af),A)$ (the algebra
of locally constant functions $\G(\Af)\fl A$ with compact support) be the
full Hecke algebra, still with the product given by convolution. Then we get
an action of $\Hecke(A)$ on the limit $\limik\ \Ho^*(S^\K)$. So the
$A$-module $\limik\ \Ho^*(S^\K)$ admits an action of the group $\G(\Af)$.

\subsection{Canonical models}
\label{subsection:SV3}

Another feature of Shimura varieties is that they have so-called
\emph{canonical
models}. That is, they are canonically defined over a number field $E$, called
the \emph{reflex field},
that depends only on $\G$ and the morphism $h:\C^\times
\fl\G(\R)$ (in particular, it does not depend on the open compact subgroup
$\K$ of $\G(\Af)$). We will use the same notation $S^\K$
for the model over $E$.
Here ``canonically'' means in particular that the action
of $\G(\Af)$ on the projective system $(S^\K)_\K$ is defined over $E$.
The theory of canonical models was begun by Shimura, and then continued by
Deligne, Borovoi, Milne and Moonen (cf \cite{D}, \cite{D2}, 
\cite{Bor}, \cite{Mi1},
\cite{Mi2}, \cite{Mo}).

So, if the cohomology theory $\Ho^*$ happens to make sense for varieties over
$E$ (for example, it could be $\ell$-adic {\'e}tale cohomology, with or
without supports), then the limit $\limik\ \Ho^*(S^\K)$ admits commuting
actions of $\G(\Af)$ and of $\Gal(\overline{E}/E)$. Another way to look at
this is to say that the cohomology group at finite level, $\Ho^*(S^\K)$,
admits commuting actions of $\Hecke_\K(A)$ and of $\Gal(\overline{E}/E)$.

The goal is now to understand the decomposition of those cohomology groups
as representations of $\G(\Af)\times\Gal(\overline{E}/E)$ (or of
$\Hecke_\K(A)\times\Gal(\overline{E}/E)$).

\subsection{Compactifications and the choice of cohomology theory}

If the Shimura varieties $S^\K$ are projective, which happens if and only if
the group $\G$ is anisotropic over $\Q$, then the most natural choice
of cohomology theory is simply the {\'e}tale cohomology of $S^\K$. There
is still the question of the coefficient group $A$. While the study of
cohomology with torsion or integral coefficients is also interesting, very
little is known about it at this point, so we will restrict ourselves to the
case $A=\overline{\Q}_\ell$, where $\ell$ is some prime number.

Things get a little more complicated when the $S^\K$ are not projective, and
this is the case we are most interested in here. We
can still use ordinary {\'e}tale cohomology or {\'e}tale cohomology with
compact support, but it becomes much harder to study (among other things,
because we do not have Poincar{\'e} duality or the fact that the cohomology
is pure - in Deligne's sense - any more). Nonetheless, it is still an
interesting problem.

Another solution is to use a cohomology theory
on a compactification of $S^\K$.
The author of this article knows of two compactifications of $S^\K$ as an
algebraic variety over $E$ (there are many, many compactifications of
$S^\K(\C)$ as a topological space, see for example the book \cite{BJ}
of Borel and Ji) :
\begin{itemize}
\item[(1)] The \emph{toroidal compactifications}. They are a family of
compactifications of $S^\K$, depending on some combinatorial data (that
depends on $\K$); they can be chosen to be very nice (i.e. projective
smooth and with a boundary that is a divisor with normal crossings).
\item[(2)] The \emph{Baily-Borel (or minimal Satake, or Satake-Baily-Borel)
compactification} $\overline{S}^\K$. It is a canonical compactification
of $S^\K$, and is a projective normal variety over $E$, but it is very
singular in general.

\end{itemize}

See the book \cite{AMRT} by Ash, Mumford, Rapoport and Tai
for the construction of the toroidal compactifications over $\C$, the
article \cite{BB} of Baily and Borel for the construction of the Baily-Borel
compactification
over $\C$, and Pink's dissertation \cite{P1} for the models over $E$ of the
compactifications.

The problem of using a cohomology theory on a toroidal compactification is
that the toroidal compactifications are not canonical, so it is not easy to
make the Hecke operators act on their cohomology. On the other hand, while
the Baily-Borel compactification is canonical (so the Hecke operators
extend to it), it is singular, so its cohomology does not behave well in
general. One solution is to use the intersection cohomology (or homology)
of the Baily-Borel
compactification. In the next section, we say a little more about intersection
homology, and explain why it might be a good choice.

\section{Intersection homology and $L^2$ cohomology}
\label{section:IH}

\subsection{Intersection homology}
\label{subsectuon:IH1}

Intersection homology was invented by Goresky and MacPherson to study the
topology of singular spaces (cf \cite{GM1}, \cite{GM2}). Let $X$ be a complex
algebraic (or analytic) variety of pure dimension $n$, possibly singular.
Then the singular homology groups of $X$ (say with coefficients in $\Q$) do
not satisfy Poincar{\'e} duality if $X$ is not smooth. To fix this, Goresky and
MacPherson modify the definition of singular homology in the following way.
First, note that $X$ admits a Whitney stratification, that is, a locally finite
decomposition into disjoint connected smooth subvarieties $(S_i)_{i\in I}$
satisfying the Whitney condition (cf \cite{GM1} 5.3). For every $i\in I$,
let $c_i=n-\dim(S_i)$ be the (complex) codimension of $S_i$.
Let $(C_k(X))_{k\in\Z}$ be the complex of simplicial chains on $X$ with
coefficients in a commutative ring $A$. The
\emph{complex of intersection chains} $(IC_k(X))_{k\in\Z}$ is the subcomplex of
$(C_k(X))_{k\in Z}$ consisting of chains $c\in C_k(X)$ satisfying the
allowability condition : For every $i\in I$, the real dimension of
$c\cap S_i$ is less than $k-c_i$, and the real dimension of $\partial c\cap
S_i$ is less than $k-1-c_i$.
The \emph{intersection homology groups}
$\IH_k(X)$ of $X$ are the homology groups of $(IC_k(X))_{k\in\Z}$.
(Note that this is the definition of middle-perversity intersection homology.
We can get other interesting intersection homology groups of $X$ by playing
with the bounds in the definition of intersection chains, but they will not
satisfy Poincar{\'e} duality.)

Intersection homology groups satisfy many of the properties of ordinary
singular
homology groups $\Ho_k(X)$ on smooth varieties.
Here are a few of these properties :
\begin{itemize}
\item[$\bullet$] They depend only on $X$, and not on the stratification
$(S_i)_{i\in I}$.
\item[$\bullet$] If $X$ is smooth, then $\IH_k(X)=\Ho_k(X)$.
\item[$\bullet$] If $X$ is compact, then the $\IH_k(X)$ are finitely generated.
\item[$\bullet$] If the coefficients $A$ are a field, the intersection
homology groups satisfy the
K{\"u}nneth theorem.
\item[$\bullet$] If $U\subset X$ is open, then there are relative intersection
homology groups $\IH_k(X,U)$ and an excision long exact sequence.
\item[$\bullet$] It is possible to define an intersection product on
intersection homology, and,
if $X$ is compact and $A$ is a field, this will induce a nondegenerate linear
pairing
\[\IH_k(X)\times \IH_{2n-k}(X)\fl A.\]
(I.e., there is a Poincar{\'e} duality theorem for intersection homology.)
\item[$\bullet$] Intersection homology
satisfies the Lefschetz hyperplane theorem and the
hard Lefschetz theorem (if $A$ is a field for hard Lefschetz).

\end{itemize}

Note however that the intersection homology groups are not homotopy invariants
(though they are functorial for certain maps of varieties, called placid maps).

\subsection{$L^2$ cohomology of Shimura varieties and intersection homology}
\label{subsection:IH2}

Consider again a Shimura variety $S^\K(\C)$ as in section
\ref{section:SV} (or rather,
the complex manifold of its complex points). For every $k\geq 0$, we write
$\Omega_{(2)}^k(S^\K(\C))$ for the space of smooth forms $\omega$ on
$S^\K(\C)$ such that $\omega$ and $d\omega$ are $L^2$.
The $L^2$ cohomology groups $\Ho_{(2)}^*(S^\K(\C))$
of $S^\K(\C)$ are the cohomology groups of the complex $\Omega_{(2)}^*$.
These groups are known to be finite-dimensional and to satisfy Poincar{\'e}
duality, and in fact we have the following theorem (remember that
$\overline{S}^\K$ is the Baily-Borel compactification of $S^\K$) :

\begin{theorem} There are isomorphisms
\[\Ho^k_{(2)}(S^\K(\C))\simeq\IH_{2d-k}(\overline{S}^\K(\C),\R),\]
where $d=\dim(S^\K)$.
Moreover, these isomorphisms are equivariant under the action
of $\Hecke_\K(\R)$.
(The Hecke algebra acts on intersection homology because the Hecke
correspondences extend to the Baily-Borel compactifications and are still
finite, hence placid.)

\end{theorem}

This was conjectured by Zucker in \cite{Z}, and then proved (independently)
by Looijenga (\cite{Lo}), Saper-Stern (\cite{SS}) and Looijenga-Rapoport
(\cite{LoR}).

So now we have some things in favour of intersection homology of the
Baily-Borel compactification : it satisfies Poincar{\'e} duality and is
isomorphic to a natural invariant of the Shimura variety. We will now see
another reason why $L^2$ cohomology of Shimura varieties (hence, intersection
homology of their Baily-Borel compactification) is easier to study than
ordinary cohomology : it is closely related to automorphic
representations of the group $\G$. (Ordinary cohomology of Shimura varieties,
or cohomology with compact support, is also related to automorphic
representations, but in a much more complicated way, see the article
\cite{F} of Franke.)

\subsection{$L^2$ cohomology of Shimura varieties and discrete automorphic
representations}
\label{subsection:IH3}

For an introduction to automorphic forms, we refer to the article \cite{BJa} of
Borel and Jacquet and the article \cite{PS} of Piatetski-Shapiro.
Let $\Ade=\Af\times\R$ be the ring of adeles of $\Q$.
Very roughly, an \emph{automorphic form} on $\G$ is a smooth function
$f:\G(\Ade)\fl\C$, left invariant under $\G(\Q)$, right invariant under
some open compact subgroup of $\G(\Af)$,
$\K_\infty$-finite on the
right (i.e., such that the right translates of $f$ by elements of $\K_\infty$
generate a finite dimensional vector space; remember that $\K_\infty$ is a
maximal compact subgroup of $\G(\R)$) and satisfying certain growth conditions.
The group $\G(\Ade)$ acts on the space of automorphic forms by right
translations on the argument. Actually, we are cheating a bit here. The group
$\G(\Af)$ does act that way, but $\G(\R)$ does not; the space of automorphic
forms is really a Harish-Chandra $(\ggoth,\K_\infty)$-module, where
$\ggoth$ is the Lie algebra of $\G(\C)$.
An \emph{automorphic representation} of $\G(\Ade)$ (or, really, $\G(\Af)\times
(\ggoth,\K_\infty)$) is an irreducible representation that appears in the
space of automorphic forms as an irreducible subquotient.

Note that there is also a classical point of view on automorphic forms, where
they are seen as smooth functions on $\G(\R)$, left invariant by some
arithmetic subgroup of $\G(\Q)$, $\K_\infty$-finite on the right and
satisfying a growth condition. From that point of view, it may be easier to
see that automorphic forms generalize classical modular forms (for
modular forms, the group $\G$ is $\GL_2$). The two points of view are closely
related, cf. \cite{BJa} 4.3 (in much
the same way that the classical and adelic
points of view on Shimura varieties are related). In this article, we adopt
the adelic point of view, because it makes it easier to see the action of
Hecke operators.

Actually, as we are interested only in discrete automorphic representations
(see below for a definition), we can see automorphic forms as $L^2$ functions
on $\G(\Q)\sous\G(\Ade)$. We follow Arthur's presentation in
\cite{A-L2}.
First, a word of warning : the quotient
$\G(\Q)\sous\G(\Ade)$ does not have finite volume. This is due to the
presence of factors isomorphic to $\R_{>0}$ in the center of $\G(\R)$.
As in \ref{subsection:SV1}, let $\Ar_\infty=\A(\R)^0$, where $\A$ is the
maximal $\R$-split torus in the center of $\G$. Then $\G(\Q)\sous\G(\Ade)/
\Ar_\infty$ does have finite volume, and we will consider $L^2$ functions
on this quotient, instead of $\G(\Q)\sous\G(\Ade)$. 

So let $\xi:\Ar_\infty\fl\C^\times$ be a character (not necessarily unitary).
Then $\xi$ extends to a character $\G(\Ade)\fl\C^\times$, that we will still
denote by $\xi$ (cf. I.3 of Arthur's introduction to the trace formula,
\cite{A-intro}).
Let $L^2(\G(\Q)\sous\G(\Ade),\xi)$ be the space of
measurable functions $f:\G(\Q)\sous\G(\Ade)\fl\C$ such that :
\begin{itemize}
\item[(1)] for every $z\in\Ar_\infty$ and $g\in\G(\Ade)$, $f(zg)=\xi(z)f(g)$;
\item[(2)] the function $\xi^{-1}f$ is square-integrable on $\G(\Q)\sous
\G(\Ade)/\Ar_\infty$.

\end{itemize}

Then the group $\G(\Ade)$ acts on $L^2(\G(\Q)\sous\G(\Ade),\xi)$ by right
translations on the argument. By definition, a \emph{discrete automorphic
representation} of $\G$ is an irreducible representation of $\G(\Ade)$ that
appears as a direct summand in $L^2(\G(\Q)\sous\G(\Ade),\xi)$. It is known
that the multiplicity of a discrete automorphic representation $\pi$ in
$L^2(\G(\Q)\sous\G(\Ade),\xi)$ is always finite; we denote it by $m(\pi)$.
We also denote by $\Pi_{disc}(\G,\xi)$ the set of discrete automorphic
representations on which $\Ar_\infty$ acts by $\xi$.
For the fact that discrete automorphic representations are indeed automorphic
representations in the previous sense, see \cite{BJa} 4.6. (The attentive reader
will have noted that automorphic representations are not actual representations
of $\G(\Ade)$ - because $\G(\R)$ does not act on them - while discrete
automorphic representations are. How to make sense of our statement that
discrete automorphic representations are automorphic is also explained in
\cite{BJa} 4.6.)

Now, given the definition of discrete automorphic representations and the
fact that $S^\K(\C)=\G(\Q)\sous\G(\Ade)/(\Ar_\infty\K_\infty\times\K)$, it
is not too surprising that the $L^2$ cohomology of the Shimura variety
$S^\K(\C)$ should be related to discrete automorphic representations.
Here is the precise relation :

\begin{theorem}(Borel-Casselman, cf. \cite{BC} theorem 4.5)
Let $\K$ be an open compact subgroup of $\G(\Af)$. Then there is a
$\Hecke_\K(\C)$-equivariant isomorphism
\[\Ho^*_{(2)}(S^\K(\C))\otimes_\R\C\simeq\bigoplus_{\pi\in\Pi_{disc}(\G,1)}
\Ho^*(\ggoth,\Ar_\infty\K_\infty;\pi_\infty)^{m(\pi)}\otimes\pi_f^\K.\]

\end{theorem}
(This is often called Matsushima's formula when $S^\K(\C)$ is compact.)

We need to explain the notation. First, the ``$1$'' in
$\Pi_{disc}(\G,1)$ stands for the trivial character of $\Ar_\infty$.
(We have chosen to work with the constant sheaf on $S^\K$, in order to
simplify the notation. In general, for a non-trivial coefficient system on
$S^\K(\C)$, other characters of $\Ar_\infty$ would appear.)
Let $\pi\in\Pi_{disc}(\G,1)$. Then $\pi$ is an irreducible representation
of $\G(\Ade)=\G(\R)\times\G(\Af)$ so it decomposes as a tensor product
$\pi_\infty\otimes\pi_f$, where $\pi_\infty$ (resp. $\pi_f$) is an irreducible
representation of $\G(\R)$ (resp. $\G(\Af)$). We denote by $\pi_f^\K$ the space
of $\K$-invariant vectors in the space of $\pi_f$; it carries an action of
the Hecke algebra $\Hecke_\K(\C)$. Finally, $\Ho^*(\ggoth,\Ar_\infty\K_\infty;
\pi_\infty)$, the $(\ggoth,\Ar_\infty\K_\infty)$-cohomology of $\pi_\infty$
(where $\ggoth$ is as before the Lie algebra of $\G(\C)$), is defined in
chapter I of the book \cite{BW} by Borel and Wallach.

This gives another reason to study the intersection homology of the
Baily-Borel compactifications of Shimura varieties : it will give a lot
of information about discrete automorphic representations of $\G$.
(Even if only about the ones whose infinite part has nontrivial
$(\ggoth,\Ar_\infty\K_\infty)$-cohomology, and that is a pretty strong
condition.)

Note that there is an issue we have been avoiding until now. Namely,
in \ref{subsection:SV3}, we wanted the cohomology theory on the Shimura
variety
to also have an action of $\Gal(\overline{E}/E)$, where $E$ is the reflex
field (i.e., the field over which the varieties $S^\K$ have canonical
models). It is not clear how to endow the $L^2$ cohomology of $S^\K(\C)$
with such an action. As we will see in the next section, this will come
from the isomorphism of $\Ho^*_{(2)}(S^\K(\C))$ with the intersection
homology of $\overline{S}^\K(\C)$ and from the sheaf-theoretic interpretation
of intersection homology (because this interpretation will also make sense
in an {\'e}tale $\ell$-adic setting).

\section{Intersection (co)homology and perverse sheaves}
\label{section:PS}

We use again the notation of section \ref{section:IH}.

\subsection{The sheaf-theoretic point of view on intersection homology}
\label{subsection:PS1}

Intersection homology of $X$ also has a sheaf-theoretical interpretation.
(At this point, we follow Goresky and MacPherson and shift from the homological
to the cohomological numbering convention.) For every open $U$ in $X$,
let $\IC^k(U)$ be the group of $(2n-k)$-dimensional intersection chains
on $U$ with closed support. If $U'\subset U$, then we have a map
$\IC^k(U)\fl \IC^k(U')$ given by restriction of chains. In this way, we get
a sheaf $\IC^k$ on $X$. Moreover, the boundary maps of the complex of
intersection chains give maps of sheaves $\delta: \IC^k\fl \IC^{k+1}$ such
that $\delta\circ\delta=0$, so the $\IC^k$
form a complex of sheaves $\IC^*$ on $X$. This is the \emph{intersection
complex} of
$X$. Its cohomology with compact support gives back the intersection homology
groups of $X$ :
\[\Ho^k_c(X,\IC^*(X))=\IH_{2n-k}(X).\]
Its cohomology groups $\IH^k(X):=\Ho^k(X,\IC^*(X))$ are (by definition) the
\emph{intersection cohomology groups} of $X$.

\subsection{Perverse sheaves}
\label{subsection:PS2}

This point of view has been extended and generalized by the invention of
perverse sheaves. The author's favourite reference for perverse sheaves is
the book by Beilinson, Bernstein and Deligne (\cite{BBD}). 

To simplify, assume
that the ring of coefficients $A$ is a field. Let $D(X)$ be the derived
category of the category of sheaves on $X$. This category is obtained from
the category of complexes of sheaves on $X$ by introducing formal inverses
of all the quasi-isomorphisms, i.e. of all the morphisms of complexes that
induce isomorphisms on the cohomology sheaves.
(This is a categorical analogue of a ring
localization.) Note that the objects of $D(X)$ are still the complexes of
sheaves, we just added more morphisms.
The homological functors on the category of complexes of
sheaves (such as the various cohomology functors and
the $Ext$ and $Tor$ functors) give functors on $D(X)$, and a morphism
in $D(X)$ is an isomorphism if and only if it is an isomorphism on the
cohomology sheaves.

This category $D(X)$ is still a little big, and we will work with the full
subcategory $D^b_c(X)$ of bounded constructible complexes. If $C^*$ is
a complex of sheaves, we will denote its cohomology sheaves by
$\Ho^kC^*$. Then $C^*$ is called \emph{bounded}
if $\Ho^kC^*=0$ for $k<<0$ and $k>>0$.
It is called \emph{constructible} if its cohomology sheaves $\Ho^kC^*$ are
constructible, that is, if, for every $k\in\Z$, there exists a
stratification $(S_i)_{i\in I}$ of $X$ (by smooth subvarieties) such that
$\Ho^kC^*_{|S_i}$ is locally constant and finitely generated for every $i$.

For every point $x$ of $X$, we denote by $i_x$ the inclusion of $x$ in $X$.

\begin{definition} A complex of sheaves $C^*$ in $D^b_c(X)$ is called
a \emph{perverse sheaf} if it satisfies the following support and
cosuport conditions :
\begin{itemize}
\item[(1)] Support : for every $k\in\Z$,
\[\dim_\C\{x\in X|\Ho^k(i_x^*C^*)\not=0\}\leq -k.\]
\item[(2)] Cosupport : for every $k\in\Z$,
\[\dim_\C\{x\in X|\Ho^k(i_x^!C^*)\not=0\}\leq k.\]

\end{itemize}

We denote by $P(X)$ the category of perverse sheaves on $X$.

\end{definition}

\begin{remark} Let $x\in X$. There is another way to look at the groups
$i_x^*\Ho^kC^*$ and $i_x^!\Ho^kC^*$. Choose an (algebraic or analytic) embedding
of a neighbourhood of $x$ into an affine space $\C^p$, and let $B_x$ denote the
intersectioon of this neighbourhood and of a small enough open ball in
$\C^p$ centered at $x$. Then
\[\Ho^k(i_x^*C^*)=\Ho^k(B_x,C^*)\]
\[\Ho^k(i_x^!C^*)=\Ho^k_c(B_x,C^*).\]

\end{remark}

\begin{remark} As before, we are only considering one perversity, the
middle (or self-dual) perversity. For other perversities (and much more), see
\cite{BBD}.

\end{remark}

Note that perverse sheaves are not sheaves but complexes of sheaves. However,
the category of perverse sheaves satisfies many properties that we expect
from a category of sheaves, and that are not true for $D^b_c(X)$ (or $D(X)$).
For example, $P(X)$ is an abelian category, and it is possible to glue
morphisms of perverse sheaves (more precisely, categories of perverse sheaves
form a stack, say on the open subsets of $X$, cf. \cite{BBD} 2.1.23).

\subsection{Intermediate extensions and the intersection complex}
\label{subsection:PS3}

Now we explain the relationship with the intersection complex. First, the
intersection complex is a perverse sheaf on $X$ once we put it in the right
degree. In fact :

\begin{prop} The intersection complex $\IC^*(X)$ is an object of $D^b_c(X)$
(i.e., it is a bounded complex with constructible cohomology sheaves), and :\begin{itemize}
\item[(1)] For every $k\not=0$,
\[\dim_\C\{x\in X|\Ho^k(i_x^*\IC^*(X))\not=0\}< n-k.\]
\item[(2)] For every $k\not=2n$,
\[\dim_\C\{x\in X|\Ho^k(i_x^!\IC^*(X))\not=0\}<k-n.\]
\item[(3)] If $U$ is a smooth open dense subset of $X$, then $\IC^*(X)_{
|U}$ is quasi-isomorphic (i.e., isomorphic in $D^b_c(X)$) to the constant sheaf
on $U$.

\end{itemize}

Moreover, the intersection complex is uniquely characterized by these properties
(up to unique isomorphism in $D^b_c(X)$).

\end{prop}

In particular, $\IC^*(X)[n]$ (that is, the intersection complex put in degree
$-n$) is a perverse sheaf on $X$.

Even better, it turns out that every perverse sheaf on $X$ is, in some sense,
built from intersection complexes on closed subvarieties of $X$. Let us be
more precise. Let $j:X\fl Y$ be a locally closed immersion. Then there is
a functor $j_{!*}:P(X)\fl P(Y)$, called the \emph{intermediate extension
functor}, such that, for every perverse sheaf $K$ on $X$, the perverse sheaf
$j_{!*}K$ on $Y$ is uniquely (up to unique quasi-isomorphism) characterized
by the following conditions :
\begin{itemize}
\item[(1)] For every $k\in\Z$,
\[\dim_\C\{x\in Y-X|\Ho^k(i_x^*j_{!*}K))\not=0\}<-k.\]
\item[(2)] For every $k\in\Z$,
\[\dim_\C\{x\in Y-X|\Ho^k(i_x^!j_{!*}K)\not=0\}< k.\]
\item[(3)] $j^*j_{!*}K=K$.

\end{itemize}

\begin{remark} Let us explain briefly the name ``intermediate extension''.
Although it is not clear from the way we defined perverse sheaves, there are
``perverse cohomology'' functors $\Hp^k:D^b_c(X)\fl P(X)$. In fact, it even
turns out that $D^b_c(X)$ is equivalent to the derived category of the abelian
category of perverse sheaves (this is a result of Beilinson, cf. \cite{Be}).
We can use these cohomology functors to define perverse extension
functors ${}^pj_!$ and ${}^pj_*$ from $P(X)$ to $P(Y)$. (For example,
${}^pj_!=\Hp^0j_!$, where $j_!:D^b_c(X)\fl D^b_c(Y)$ is the ``extension
by zero'' functor between the derived categories; likewise for ${}^pj_*$).
It turns out that, from the perverse point of view, the functor
$j_!:D^b_c(Y)\fl D^b_c(X)$ is right exact and the functor $j_*:D^b_c(Y)\fl
D^b_c(X)$ is left exact (that, if $K$ is perverse on $X$, $\Hp^kj_!K=0$
for $k>0$ and $\Hp^kj_*K=0$ for $k<0$). So the morphism of functors $j_!\fl
j_*$ induces a morphism of functors ${}^pj_!\fl{}^pj_*$. For every perverse
sheaf $K$ on $X$, we have :
\[j_{!*}K=\Ima({}^pj_!K\fl{}^pj_*K).\]

\end{remark}

Now we come back to the description of the category of perverse sheaves on
$X$. Let $F$ be a smooth connected
locally closed subvariety of $X$, and denote by
$i_F$ its inclusion in $X$. If $\F$ is a locally constant sheaf on $F$, then
it is easy to see that $\F[\dim F]$ is a perverse sheaf on $F$; so
$i_{F!*}\F[\dim F]$ is a perverse sheaf on $X$ (it has support in
$\overline{F}$, where $\overline{F}$ is the closure of $F$ in $X$).
If the locally constant
sheaf $\F$ happens to be irreducible, then this perverse sheaf is a simple
object in $P(X)$. In fact : 

\begin{theorem} The abelian category $P(X)$ is artinian and noetherian (i.e.,
every object has finite length), and its simple objects are all of the form
$i_{F!*}\F[\dim F]$, where $F$ is as above and $\F$ is an irreduible locally
constant sheaf on $F$.

\end{theorem}

Finally, here is the relationship with the intersection complex. Let
$i_F:F\fl X$ be as above. Then, if $\F$ is the constant sheaf on $F$,
the restriction to $\overline{F}$ of the perverse sheaf $i_{F!*}\F[\dim F]$
is isomorphic to $\IC^*(\overline{F})[\dim F]$. In fact, we could define the
intersection complex on a (possibly singular) variety $Y$ with coefficients
in some locally constant sheaf on the smooth locus of $Y$, and then the simple
objects in $P(X)$ would all be intersection complexes on closed
subvarieties of $X$.

\subsection{$\ell$-adic perverse sheaves}

Now we come at last to the point of this section (to make the Galois
groups $\Gal(\overline{E}/E)$ act on the intersection (co)homology of
$\overline{S}^\K(\C)$).

Note that the definitions of the category of perverse sheaves and of the
intermediate extension in \ref{subsection:PS2} and \ref{subsection:PS3} would
work just as well in a category of {\'e}tale $\ell$-adic sheaves. So now
we take for $X$ a quasi-separated scheme of finite type over a field
$k$, we fix a prime number $\ell$ invertible in $k$
and we consider the category $D^b_c(X,\overline{\Q}_\ell)$ of bounded
$\ell$-adic complexes on $X$.
(To avoid a headache, we will take $k$ to be algebraically closed or finite,
so the simple construction of \cite{BBD} 2.2.14 applies.)
Then we can define an abelian subcategory of perverse sheaves $P(X)$ in
$D^b_c(X,\overline{\Q}_\ell)$ and intermediate extension functors
$j_{!*}:P(X)\fl P(Y)$ as before (see \cite{BBD} 2.2). In particular, we can make
the following definition :

\begin{definition} Suppose that $X$ is purely of dimension $n$, and let
$j:U\fl X$ be the inclusion of the smooth locus of $X$ in $X$. Then the
($\ell$-adic) intersection complex of $X$ is
\[\IC^*(X)=(j_{!*}\overline{\Q}_{\ell,U}[n])[-n],\]
where $\overline{\Q}_{\ell,U}$ is the constant sheaf $\overline{\Q}_\ell$ on
$U$. The $\ell$-adic intersection cohomology $\IH^*(X,\overline{\Q}_\ell)$
of $X$ is the cohomology of $\IC^*(X)$.

\end{definition}

\subsection{Application to Shimura varieties}

We know that the Shimura variety $S^\K$ and its Baily-Borel compactification
$\overline{S}^\K$ are defined over the number field $E$. So we can form the
$\ell$-adic intersection cohomology groups
$\IH^*(\overline{S}^\K_{\overline{E}},
\overline{\Q}_\ell)$. They admit an action of $\Gal(\overline{E}/E)$. Moreover,
if we choose a field isomorphism $\overline{\Q}_\ell\simeq\C$, then the
comparison theorems between the {\'e}tale topology and the classical
topology will give an isomorphism $\IH^*(\overline{S}^\K_{\overline{E}},
\overline{\Q}_\ell)\simeq\IH^*(\overline{S}^\K(\C),\C)$ (cf. chapter 6 of
\cite{BBD}).

The isomorphism of \ref{subsection:IH2}
between intersection homology of $\overline{S}^\K(\C)$
and $L^2$ cohomology of $S^\K(\C)$, as well as the duality between
intersection homology and intersection cohomology (cf. \ref{subsection:PS1}),
thus give an isomorphism
\[\IH^*(\overline{S}^\K_{\overline{E}},\overline{\Q}_\ell)\simeq
\Ho^*_{(2)}(S^\K(\C))\otimes\C,\]
and this isomorphism is equivariant under the action of $\Hecke_\K(\C)$.
We know what $L^2$ cohomology looks like as a representation of $\Hecke_\K(\C)$,
thanks to the theorem of Borel and Casselman (cf. \ref{subsection:IH3}).

Using this theorem and his own trace
invariant formula, Arthur has given a formula
for the trace of a Hecke operator on $\Ho^*_{(2)}(S^\K(\C))\otimes\C$ (cf.
\cite{A-L2}). This formula involves global
volume terms, discrete series characters
on $\G(\R)$ and orbital integrals on $\G(\Af)$.

The
problem now is to understand the action of the Galois group $\Gal(\overline{E}/
E)$. We have a very precise conjectural description of the intersection
cohomology of $\overline{S}^\K$ as a $\Hecke_\K(\C)\times\Gal(\overline
{E}/E)$-module, see for example the articles \cite{K-SVLR} of Kottwitz
and \cite{BR} of Blasius and Rogawski.

In the next sections, we will explain a strategy to understand how at least
part of the Galois group $\Gal(\overline{E}/E)$ acts.

\section{Counting points on Shimura varieties}
\label{section:CP}

We want to understand the action of the Galois
group $\Gal(\overline{E}/E)$ on the intersection cohomology groups
$\IH^*_\K:=\IH^*(\overline{S}^\K_{\overline{E}},\overline{\Q}_\ell)$. It is
conjectured that this action is unramified almost everywhere. Thus, by the
Chebotarev density theorem,
it is theoretically enough to understand the action of
the Frobenius automorphisms at the places of $E$ where the action is
unramified, and one way to do this is to calculate the trace of the
powers of the Frobenius
automorphisms at these places. However, for some purposes, it is necessary to
look at the action of the decomposition groups at other places. This is part
of the theory of bad reduction of Shimura varieties, and we will not talk
about this here, nor will we attempt to give comprehensive references to it.
(Let us just point to the book \cite{HT} of Harris and Taylor.)

In general, intersection cohomology can be very hard to calculate. First
we will look at simpler objects, the cohomology groups with compact support
$\Ho^*_{c,\K}:=\Ho_c^*(S^\K_{\overline{E}},\overline{\Q}_\ell)$.
Assume that the Shimura varieties and their compactifications (the Baily-Borel
compactifications and the toroidal compactifications) have ``good'' models
over an open subset $U$ of $\Spec\Of_E$,
and write $\Scal^\K$ for the model of $S^\K$. (It is much easier to imagine what
a ``good'' model should be than to write down a precise definition. An
attempt has been made in \cite{Mo2} 1.3, but it is by no means optimal.)
Then, by the specialization theorem (SGA 4 III Expos{\'e} XVI 2.1), and also
by Poincar{\'e} duality (cf. SGA 4 III Expos{\'e} XVIII), for every finite
place $\pp$ of $E$ such that $\pp\in U$ and $\pp\not|\ell$,
there is a $\Gal(\overline{E}_\pp/E_\pp)$-equivariant isomorphism
\[\Ho^*_{c,\K}=\Ho_c^*(S^\K_{\overline{E}},\overline{\Q}_\ell)\simeq
\Ho_c^*(\Scal^\K_{\overline{\Fi}_\pp},\overline{\Q}_\ell),\]
where $\Fi_\pp$ is the residue field of $\Of_E$ at $\pp$.
In particular, the $\Gal(\overline{E}/E)$-representation $\Ho^*_{c,\K}$ is
unramified at $\pp$.

Now, by Grothendieck's fixed point formula (SGA 4 1/2 Rapport), calculating
the trace of powers of the Frobenius automorphism on $\Ho_c^*(\Scal^\K_
{\overline{\Fi}_\pp},\overline{\Q}_\ell)$ is the same as counting the points
of $\Scal^\K$ over finite extensions of $\Fi_\pp$.

Langlands has given a conjectural formula for this number of points, cf.
\cite{La} and \cite{K-SVLR}. Ihara had earlier made and proved a similar
conjecture for Shimura varieties of dimension $1$. 
Although this conjecture is not known in general, it is easier to study for
a special class of Shimura varieties, the so-called PEL Shimura varieties.
These are Shimura varieties that can be seen as moduli spaces of abelian
with certain supplementary structures (P : polarizations, E : endomorphisms,
i.e. complex multiplication by certain CM number fields, and L : level
structures). For PEL Shimura varieties of types $A$ and $C$ (i.e., such that
the group $\G$ is of type $A$ or $C$), Langlands's conjecture had been proved
by Kottwitz in \cite{K-PSSV}. Note that all the examples we gave
in \ref{subsection:SV1} are of this type.
Conveniently enough, the modular interpretation of PEL Shimura varieties also
gives a model of the Shimura variety over an explicit
open subset of $\Spec\Of_E$.

In fact, Kottwitz has done more than counting
points; he has also counted the points that are fixed by the composition of
a power of the Frobenius automorphism and of a Hecke
correspondence (with a condition of triviality at $\pp$). So, using Deligne's
conjecture instead of Grothendieck's fixed point formula, we can use
Kottwitz's result to understand the commutating actions of $\Gal(\overline{E}/
E)$ and of $\Hecke_\K(\overline{\Q}_\ell)$ on $\Ho^*_{c,\K}$. (Deligne's
conjecture gives a simple formula for the local terms in the Lefschetz
fixed formula if we twist the correspondence by a high power of the Frobenius.
It is now a theorem and has been proved independently by Fujiwara in
\cite{Fu} and Varshavsky in \cite{V}. In the case of Shimura varieties, it
also follows from an earlier result of Pink in \cite{P3}.)

Using his counting result, Kottwitz has proved the conjectural description
of $\IH^*_{\K}$ for some simple Shimura varieties (cf. \cite{K-LAR}).
Here ``simple'' means that the Shimura varieties are compact (so intersection
cohomology is cohomology with compact support) and that the phenomenon called
``endoscopy'' (about which we are trying to say as little as possible) does
not appear. 

One reason to avoid endoscopic complications was that a very important
and necessary result when dealing with endoscopy,
the so-called ``fundamental lemma'', was not available at the time.
It now is, thanks to the combined efforts of many people, among
which Kottwitz (\cite{K-BC}), Clozel (\cite{C-BC}), Labesse (\cite{La-BC},
\cite{CL}), Hales (\cite{Ha}), Laumon, Ngo (\cite{LN}, \cite{Ng}),
and Waldspurger (\cite{Wa1}, \cite{Wa2}, \cite{Wa3}). 

Assuming the fundamental lemma, 
the more general case of compact PEL Shimura varieties of type $A$ or $C$
(with endoscopy playing a role) was treated by Kottwitz in \cite{K-SVLR},
admitting Arthur's conjectures on the descripton of discrete automorphic
representations of $\G$. Actually, Kottwitz did more : he treated the
case of the (expected) contribution of $\Ho^*_{c,\K}$ to $\IH^*_\K$.
Let us say a word about Arthur's conjectures. Arthur has announced a proof
of a suitable formulation of his conjectures for classical groups (that is,
symplectic and orthogonal groups), using the stable twisted trace formula.
His proof is expected to adapt to the case of unitary groups (that is, the
groups that give PEL Shimura varieties of type $A$), but this
adaptation will likely require a lot of effort.

Let us also note that the case of compact PEL Shimura varieties of type $A$
should be explained in great detail in the book
project led by Michael Harris (\cite{BP}). 

This does not tell us what to do in the case where $S^\K$ is not projective.
First note that the modular interpretation gives us integral models
of the Shimura varieties but not of their compactifications. So this is
the first problem to solve. Fortunately, it has been solved :
See the article \cite{DR} of Deligne and Rapoport for the case of modular
curves, the book \cite{CF} by Chai and Faltings for the case of
Siegel modular varieties,
Larsen's article \cite{Lar} for the case of Picard modular
varieties, and Lan's dissertation \cite{Lan} for the general case of
PEL Shimura varieties of type $A$ or $C$.
This allows us to apply the specialization theorem to intersection cohomology.
In particular, we get the fact that the $\Gal(\overline{E}/E)$-representation
$\IH^*_{c,\K}$ is unramified almost everywhere, and, at the finite places $\pp$
where it is unramified, we can study it by considering the reduction modulo
$\pp$ of the Shimura variety and its compactifications.

Next we have to somehow describe the intersection complex.
If the group $\G$ has semi-simple $\Q$-rank $1$, so it has only one conjugacy
class of rational
parabolic subgroups, then the Baily-Borel compactification is simpler
(it only has one kind of boundary strata) and
we can obtain the intersection complex
by a simple truncation process from the direct image on $\overline{S}^\K$
of the constant sheaf on $S^\K$. The conjectural description of $\IH^*_\K$
is know for the cases $\G=\GL_2$ (see the book \cite{DK}) and the case
of Picard modular surfaces, i.e., $\G=\GU(2,1)$ (see the book \cite{LR}).
In the general case of semi-simple $\Q$-rank $1$, Rapoport has given in
\cite{Ra} a formula for the trace of a power of the Frobenius automorphism
(at almost every place) on the stalks of the intersection complex.

In the general case, the intersection complex  is obtained from the direct
image of the constant sheaf on $S^\K$ by applying several nested truncations
(cf. \cite{BBD} 2.1.11), and it is not clear how to see the action of Frobenius
on the stalks of this thing. We will describe a solution in the next section.

\section{Weighted cohomology}
\label{section:WC}

In this section, $j$ will be the inclusion of $S^\K$ in its Baily-Borel
compactification $\overline{S}^\K$, and $j_*$ will be the derived direct
image functor.
Here is the main idea : instead of seeing the intersection complex
$IC^*(\overline{S}^\K)$ as a truncation of $j_*\overline{\Q}_{\ell,S^\K}$
by the cohomology degree (on various strata of $\overline{S}^\K-S^\K$), we want
to see it as a truncation by Frobenius weights (in the sense of Deligne).
This idea goes back to the
construction by Goresky, Harder and MacPherson of the weighted cohomology
complexes in a topological setting (i.e., on a non-algebraic compactification
of the set of complex points $S^\K(\C)$). 

\subsection{The topological case}

As we have mentioned before, the manifold $S^\K(\C)$ has a lot of
non-algebraic compactifications (these compactifications are defined for a
general locally symmetric space, and not just for a Shimura variety).
The one used in the construction of
weighted cohomology is the reductive Borel-Serre compactification
$S^\K(\C)^{RBS}$ (cf. \cite{BJ} III.6 and III.10;
the reductive Borel-Serre compactification
was originally defined by Zucker in \cite{Z}, though not under that name).
The reductive Borel-Serre compactification admits a map
$\pi:S^\K(\C)^{RBS}\fl\overline{S}^\K(\C)$ that extends the identity on
$S^\K(\C)$; we also denote by $\til{j}$ the inclusion of
$S^\K(\C)$ in $S^\K(\C)^{RBS}$.

The boundary $S^\K(\C)^{RBS}-S^\K(\C)$ of $S^\K(\C)^{RBS}$ has a very
pleasant description. It is a union of strata, each of which is a locally
symmetric space for the Levi quotient of a rational parabolic subgroup of
$\G$; moreover, the closure of a stratum is its reductive Borel-Serre
compactification. (A lot more is known about the precise geometry of the
strata, see, e.g., \cite{GHM} 1D).

The weighted cohomology complexes are bounded constructible complexes
$W^\mu$ of $\C$ or $\Q$-vector spaces on $S^\K(\C)^{RBS}$ extending
the constant sheaf on $S^\K(\C)$, constructed by Goresky, Harder and
MacPherson in \cite{GHM} (they give two constructions, one for
$\C$-coefficients and one for $\Q$-coefficients, and then show that the two
constructions agree). They depend on a weight profile $\mu$ (which is a
function from the set of relative simple roots of $\G$ to $\Z+\frac{1}{2}$).
The basic idea of weighted cohomology is to consider the complex
$\til{j}_*\C$ (or $\til{j}_*\Q$) on $S^\K(\C)^{RBS}$ and to truncate it,
not by the cohomology degree as for the intersection complex, but by
the weights of certain tori. More precisely, on a strata $S$ corresponding to
a Levi subgroup $\M$, we truncate by the weights of the $\Q$-split torus
$\A_M$
in the center of $\M$ (the group $\A_M(\Q)$ acts
on $\til{j}_*\C_{|S}$ by what Goresky, Harder
and MacPherson call Looijenga Hecke correspondences). The weight profile
specifies, for every strata, which weights to keep.

Of course, it is not that simple. The complex $\til{j}_*\C$ is an object
in a derived category (which is not abelian but triangulated),
and it is not so easy to truncate objects in such a category. To get
around this problem, the authors of \cite{GHM} construct an incarnation
of $\til{j}_*\C$, that is, an explicit complex that is quasi-isomorphic
to $\til{j}*\C$ and on which the tori $\A_M(\Q)$ still act.
(In fact, they construct two incarnations, one of
$\til{j}_*\C$ and one of $\til{j}_*\Q$).

The upshot (for us) is that the functor $\pi_*:D^b_c(S^\K(\C)^{RBS})
\fl D^b_c(\overline{S}^\K(\C))$ sends two of these weighted cohomology
complexes to the intersection complex on $\overline{S}^\K(\C)$ (they are
the complexes corresponding to the lower and upper middle weight profiles).
On the other hand, the weighted cohomology complexes are canonical enough so
that the Hecke algebra acts on their cohomology, and explicit enough so that
it is possible to calculate the local terms when we apply the Lefschetz fixed
point formula to them. This is possible but by no means easy, and is the
object of the article \cite{GM3} of Goresky and MacPherson. Then, in the
paper \cite{GKM}, Goresky, Kottwitz and MacPherson show that the result of
\cite{GM3} agrees with the result of Arthur's calculation in \cite{A-L2}.

The problem, from our point of view, is that this construction is absolutely
not algebraic, so it is unclear how to use it to understand the action of
$\Gal(\overline{E}/E)$ on $\IH^*(S^\K,\overline{\Q}_\ell)$.

\begin{remark} There is another version of weighted cohomology of
locally symmetric spaces : Franke's weighted $L^2$ cohomology, defined in
\cite{F}. In his article \cite{Na}, Nair has shown that Franke's weighted
$L^2$ cohomology groups are weighted cohomology groups in the sense of
Goresky-Harder-MacPherson.

\end{remark}

\subsection{Algebraic construction of weighted cohomology}

First, the reductive Borel-Serre compactification is not an algebraic variety,
so what we are really looking for is a construction of the complexes
$\pi_*W^\mu$, directly on the Baily-Borel compactification. This looks
difficult for several reasons. The Baily-Borel compactification
is very singular, which is one
of the reasons why Goresky, Harder and MacPherson use the less singular
reductive Borel-Serre compactification in the first place. Besides, the
boundary strata in $\overline{S}^\K$ correspond to maximal rational parabolic
subgroups of $\G$, and several strata in $S^\K(\C)^{RBS}$ can be
(rather brutally) contracted to the same stratum in $\overline{S}^\K(\C)$.
It is possible to give a description of the stalks of $\pi_*W^\mu$
(see the article \cite{GHMN} of Goresky, Harder, MacPherson and Nair), but
it is a rather complicated description, much more complicated than the
simple description of the stalks of $W^\mu$.

The idea is that the action of the Looijenga Hecke
correspondences should correspond in some way to the action of the
Frobenius automorphism in an algebraic setting. 
This is actually a very natural ideal. Looijenga himself uses the fact that
the eigenspaces of the Looijenga Hecke correspondences are pure in the sense
of mixed Hodge theory (cf. \cite{Lo} 4.2), and we know that the weight
filtration of Hodge theory corresponds to the filtration by Frobenius
weights in $\ell$-adic cohomology (cf. for example \cite{BBD} 6.2.2).
So the correct algebraic
analogue of the truncations of \cite{GHM} should be a truncation by Frobenius
weights (in the sense of Deligne's \cite{D3}, see also chapter 5 of
\cite{BBD}).
As a consequence, the most natural place to define the algebraic analogues
of the weighted cohomology complexes is the reduction modulo $\pp$ of
an integral model of $\overline{S}^\K$, where $\pp$ is a finite place of
$E$ where good integral models exist. (But see the remark at the end of
this subsection.)

In fact, it turns out that we can work in a very general setting. Let
$\Fi_q$ be a finite field, and $X$ be a quasi-separated scheme of finite
type over $\Fi_q$. Then we have the category of mixed $\ell$-adic
complexes $D^b_m(X,\overline{\Q}_\ell)$ on $X$, cf. \cite{BBD} 5.1.
(Here ``mixed'' refers to the weights of the complexes, and the weights are
defined by considering the action of the Frobenius automorphisms on the
stalks of the complexes; for more details, see \cite{D3} or \cite{BBD} 5).
In particular, we get a category $P_m(X)$ of mixed $\ell$-adic
perverse sheaves on $X$ as a subcategory of $D^b_m(X,\overline{\Q}_\ell)$.
One important result of the theory is that mixed perverse sheaves admit
a canonical weight filtration. That is, if $K$ is an object in
$P_m(X)$, then it has a canonical filtration $(w_{\leq a}K)_{a\in\Z}$ such that
each $w_{\leq a}K$ is a subperverse sheaf of $K$ of weight $\leq a$ and such
that $K/w_{\leq a}K$ is of weight $>a$.

This functor $w_{\leq a}$ on mixed perverse sheaves
does not extend to $D^b_m(X,\overline{\Q}_\ell)$ in
the na{\"i}ve way; that is, the inclusion functor from the category
of mixed sheaves of weight $\leq a$ to $D^b_m(X,\overline{\Q}_\ell)$ does not
admit a right adjoint. But we can extend $w_{\leq a}$ in another way. Consider
the full subcategory $\DP^{\leq a}$ of $D^b_m(X,\overline{\Q}_\ell)$ whose
objects are the complexes $K$ such that, for every $k\in\Z$, the
$k$-th perverse cohomology sheaf $\Hp^kK$ is of weight $\leq a$. (If we wanted
to define the complexes of weight $\leq a$, we would require
$\Hp^kK$ to be of weight $\leq a+k$.) Then $\DP^{\leq a}$ is a
triangulated subcategory of $D^b_m(X,\overline{\Q}_\ell)$, and the inclusion
$\DP^{\leq a}\subset D^b_m(X,\overline{\Q}_\ell)$ does admit a right
adjoint, which we denote by $w_{\leq a}$ (because it extends the previous
$w_{\leq a}$). Likewise, we can define a full triangulated subcategory
$\DP^{\geq a}$ of $D^b_m(X,\overline{\Q}_\ell)$, whose inclusion into
$D^b_m(X,\overline{\Q}_\ell)$ admits a left adjoint $w_{\geq a}$ (extending
the functor $K\fle K/w_{\leq a-1}K$ on mixed perverse sheaves).
This is explained in section 3 of \cite{Mo1}.
Then the analogue of the theorem that $\pi_*W^\mu$ is the intersection
complex (for a well-chosen weight profile $\mu$) is the :

\begin{theorem}(\cite{Mo1} 3.1.4) Let $j:U\fl X$ a nonempty open
subset of $X$ and $K$ be a pure perverse sheaf of weight $a$ on $U$.
Then there are canonical isomorphisms :
\[j_{!*}K\simeq w_{\leq a}j_*K\simeq w_{\geq a}j_!K.\]

\end{theorem}

More generally, if we have a stratification on $X$, we can choose to truncate
by different weights on the different strata (cf. \cite{Mo1} 3.3); in this way,
we get analogues of the other weighted cohomology complexes, or rather
of their images on the Baily-Borel compactification. We also get somewhat
more explicit formulas for $w_{\leq a}$, and hence the intersection
complex (\cite{Mo1} 3.3.4 and 3.3.5), analogous to the formula
of \cite{BBD} 2.1.11, but where all the truncations by the cohomology
degree have been replaced by weight truncations. The reason this makes
such a big difference is that the weight truncation functors $w_{\leq a}$
and $w_{\geq a}$ are exact in the perverse sense. (Interestingly enough, it
turns out that, in this setting, the weighted cohomology complexes are
canonically defined and have nothing to do with Shimura varieties.
In fact, there is another application of these ideas, to Schubert varieties,
see \cite{Mo3}.)

\begin{remark} We want to make a remark about the construction of
the weighted cohomology complexes on the canonical models
$\overline{S}^\K$ (and not their reduction modulo a prime ideal).
The construction of \cite{Mo1} 3 is very formal and will apply in every
category that has a notion of weights and a weight truncation on ``perverse''
objects. For example, it should apply without any changes
to Saito's derived category of
mixed Hodge modules.
In fact, Arvind Nair has just informed the author that he has indeed been
able to construct weighted cohomology complexes in the category of mixed
Hodge modules, and to prove that the weighted cohomology complexes he obtained
on the Baily-Borel compactification of a Shimura variety
are the pushforwards of the
Goresky-Harder-MacPherson weighted cohomology complexes on the reductive
Borel-Serre compactification. As an application of this, he was able to
prove that Franke's spectral sequence (\cite{F} 7.4) is a spectral
sequence of mixed Hodge structures (for the locally symmetric spaces that are
Shimura varieties).

Now suppose that $X$ is a quasi-separated scheme of finite type over a number
field.
We can define $\ell$-adic perverse sheaves on
$X$, and we can also define a notion of weights for $\ell$-adic complexes on
$X$ (cf. Deligne's \cite{D3} 1.2.2 and
Huber's article \cite{Hu}). The problem is that mixed perverse sheaves
on $X$ do not have a weight filtration in general (because number fields have
more Galois cohomology than finite fields). To circumvent this problem,
we could try to work in the derived category of the abelian category
of mixed perverse sheaves on $X$ admitting a weight filtration. Then
it is not obvious how to construct the 4/5/6 operations on these
categories. It might be possible to copy Saito's approach in
\cite{Sa1} (where he constructs and studies the derived category of mixed
Hodge modules); see also Saito's preprint \cite{Sa2}. As far as the author
knows, this has not been worked out anywhere.

\end{remark}

\subsection{Application to the cohomology of Shimura varieties}

Once we have the interpretation of the intermediate extension functor
given in the previous subsection,
it becomes surprisingly easy to calculate the trace of Frobenius automorphisms
on the stalks of $\IC^*(\overline{S}^\K)$. We should mention that one
reason it is so easy is that one of the main ingredients, a description of
the restriction to the boundary strata of the complex $j_{*}\overline{\Q}_\ell$
(where $j$ is again the inclusion of $S^\K$ in $\overline{S}^\K$) has been
provided by Pink in \cite{P2}. And of course, the whole calculation rests
on Kottwitz's calculations for the cohomology with compact support
(in \cite{K-PSSV}). Including Hecke correspondences in the picture is just
a matter of bookkeeping, and the final result of the Lefschetz trace
formula appears in \cite{Mo2} 1.7.

This is not the end of the story. It still remains to compare the result
of the Lefschetz fixed point formula with Arthur's invariant trace formula,
in order to try to prove the result conjectured in 10.1 of Kottwitz's article
\cite{K-SVLR}. This is basically a generalization of part I of
\cite{K-SVLR} to include the non-elliptic terms. Given the work done by
Kottwitz in \cite{K-SVLR} and \cite{K-NP}, it requires no new ideas, but
still takes some effort. In the case of general unitary groups over $\Q$,
it is the main object of the book \cite{Mo2} (along with some applications).

Even then, we are not quite done. If we want to prove the conjectural
description of $\IH^*(\overline{S}^\K,\overline{\Q}_\ell)$ given
in \cite{K-SVLR} or \cite{BR}, we still need to know Arthur's conjectures.

Some applications that do not depend on Arthur's conjectures are worked out
in the book \cite{Mo2} (subsection 8.4).
They use a weak form of base change from unitary groups
to general linear groups, for the automorphic representations that appear in
the $L^2$ cohomology of Shimura varieties.
(If we knew full base change, then we would probably
also know Arthur's conjectures.)
Let us mention the two main applications :
\begin{itemize}
\item[$\bullet$] The logarithm of the $L$-function
of the intersection complex is a linear combination of logarithms of
$L$-functions of automorphic representations of general linear groups (\cite{Mo2}
corollary 8.4.5). In fact, we can even get similar formulas for the
$L$-functions of the $\Hecke_\K(\overline{\Q}_\ell)$-isotypical components
of the intersection cohomology, as in \cite{Mo2} 7.2.2.
However, the coefficients in these linear combinations are not explicit, and
in particular \cite{Mo2} does not show that they are integers.

\item[$\bullet$] We can derive some cases of the global Langlands
correspondence (cf. \cite{Mo2} 8.4.9, 8.4.10). Note however that one of the
conclusions of \cite{Mo2} is that, in the end, we do not get more Galois
representations in the cohomology of noncompact unitary varieties than we
would in the cohomology of compact unitary Shimura varieties. In particular,
the cases of the Langlands correspondence that are worked out in \cite{Mo2}
can also be obtained using compact Shimura varieties and gluing of Galois
representations (cf. the last chapters of the book project \cite{BP} or
the article \cite{Sh} of Shin; note that Shin also considers places of
bad reduction).

\end{itemize}

\end{document}